\documentclass[graybox,footinfo]{svmult}

\smartqed
\usepackage{mathptmx}       
\usepackage{helvet}         
\usepackage{courier}        
\usepackage{type1cm}        
\usepackage{graphicx}       

\usepackage{array,colortbl}
\usepackage{amsmath,amsfonts,amssymb,bm} 
\usepackage{url}

\newcommand{\Z}{\mathbb{Z}} 
\newcommand{\R}{\mathbb{R}} 
\newcommand{\N}{\mathbb{N}} 
\newcommand{\rd}{\,\mathrm{d}} 
\newcommand{\bsx}{\boldsymbol{x}}    

\newcommand{\calP}{\mathcal{P}}
\newcommand{\calF}{\mathcal{F}}
\newcommand{\calI}{\mathcal{I}}
\newcommand{\wafom}{{{\mathrm{WF}}}}
\newcommand{\Vol}{{{\mathrm{Vol}}}}
\newcommand{\Ftwo}{{\mathbb{F}_2}}
\def\Error{{{\mathrm{Error}}}}

\makeatletter 
  \providecommand*{\toclevel@author}{999}
  \providecommand*{\toclevel@title}{0}
\makeatother

\begin{document}

\title*{Walsh Figure of Merit for Digital Nets: An Easy Measure
for Higher Order Convergent QMC}
\titlerunning{WAFOM for Higher Order Convergent QMC}
\author{Makoto Matsumoto \and Ryuichi Ohori}
\institute{
Makoto Matsumoto
\at
Graduate School of Sciences,
Hiroshima University, Hiroshima 739-8526 Japan,
\email{m-mat@math.sci.hiroshima-u.ac.jp}
\and 
Ryuichi Ohori
\at
Fujitsu Laboratories Ltd.,
Kanagawa 211-8588 Japan,
\email{ohori.ryuichi@jp.fujitsu.com}
}
\maketitle


\abstract{
Fix an integer $s$. Let $f:[0,1)^s \rightarrow \R$ 
be an integrable function.
Let $P\subset [0,1]^s$ be a finite point set.
Quasi-Monte Carlo integration of $f$ by $P$ is
the average value of $f$ over $P$ that approximates
the integration of $f$ over the $s$-dimensional cube.
Koksma-Hlawka inequality tells that, by a smart choice of $P$, one may expect that
the error decreases roughly $O(N^{-1}(\log N)^s)$. 
For any $\alpha \geq 1$, J.\ Dick gave a construction of point sets 
such that for $\alpha$-smooth $f$, 
convergence rate $O(N^{-\alpha}(\log N)^{s\alpha})$ is assured.
As a coarse version of his theory, M-Saito-Matoba introduced
Walsh figure of Merit (WAFOM), which gives the convergence rate
$O(N^{-C\log N/s})$. WAFOM is efficiently computable.
By a brute-force search of low WAFOM point sets, 
we observe a convergence rate of order $N^{-\alpha}$
with $\alpha>1$, for several test integrands for $s=4$ and $8$.
}

\section{Quasi-Monte Carlo and Higher Order Convergence}\label{sec:1}
Fix an integer $s$. Let $f:[0,1)^s \rightarrow \R$ 
be an integrable function.
Our goal is to have a good approximation of the value 
$$I(f):=\int_{[0,1)^s} f(x) \D x.$$
We choose a finite point set 
$\calP \subset [0,1)^s$,
whose cardinality
is called the sample size and denoted by $N$.
The quasi-Monte Carlo (QMC) integration of $f$ 
by $\calP$ is the value
$$I(f;\calP):=\frac{1}{N} \sum_{x \in \calP} f(x),$$
i.e., the average of $f$ over the finite points $\calP$
that approximates $I(f)$.
The QMC integration error is defined by
$$\Error(f;\calP):=|I(f) - I(f;\calP)|.$$
If $\calP$ consists of $N$ independently, uniformly and randomly 
chosen points, the QMC integration is nothing but the classical 
Monte Carlo (MC) integration, where the integration error is expected
to decrease with the order of $N^{-1/2}$ when $N$ increases, 
if $f$ has a finite variance.

The main purpose of QMC integration is to choose good point sets
so that the integration error decreases faster than MC.
There are enormous studies in diverse directions,
see for examples \cite{DICK-PILL-BOOK} \cite{niederreiter:book}.

In applications, often we know little on the integrand $f$,
so we want point sets which work well for a wide class
of $f$.
An inequality
of the form
\begin{equation}\label{eq:koksma}
\Error(f;\calP)  \le V(f) D(\calP),
\end{equation}
called of Koksma-Hlawka type, is often useful.
Here, $V(f)$ is a value independent of $\calP$
which measures some kind of variance of $f$, and
$D(\calP)$ is a value independent of $f$
which measures some kind of discrepancy 
of $\calP$ from an ``ideal'' uniform distribution.
Under such an inequality,
we may prepare point sets 
with small values of $D(\calP)$, 
and use them for QMC-integration
if $V(f)$ is expected to be not too large.

In the case of the original Koksma-Hlawka inequality,
\cite[Chapters 2 and 3]{niederreiter:book},
$V(f)$ is the total variation
of $f$ in the sense of Hardy and Krause,
and $D(\calP)$ is the star discrepancy of 
the point set.
In this case the inequality is known to be sharp.
It is a conjecture that there is a constant $c_s$
depending only on $s$ such that
$D^*(\calP)>c_s (\log N)^{s-1}/N$,
and there are constructions of point sets
with 
$D^*(\calP)<C_s (\log N)^{s}/N$.
Thus, to obtain a better convergence rate,
one needs to assume some restriction on $f$.
If for a function class $\calF$, there
are $V(f)$ $(f\in \calF)$ and $D(\calP)$ with 
the inequality (\ref{eq:koksma}) with a sequence of 
point sets $\calP_1, \calP_2, \ldots$ with
$D(\calP_i)$ decreases faster than the order
$1/N_i$, then it is natural to call
the point sets as higher order QMC point sets
for the function class $\calF$.

It is known that this is possible if we assume
some smoothness on $f$. 
Dick \cite{Dick-Walsh} \cite{Dick-MCQMC} \cite{DICK-PILL-BOOK}
showed that for any positive integer $\alpha$,
there is a function class named $\alpha$-smooth 
such that
the inequality
$$
\Error(f; \calP) \leq C(\alpha,s)||f||_\alpha W_\alpha(\calP)
$$
holds,
where point sets with
$
W_\alpha(\calP)=O(N^{-\alpha}(\log N)^{s\alpha})
$ are constructible 
from $(t,m,s)$-nets 
(named higher order digital net).
The definition of $W_\alpha(\calP)$ is given later
in \S\ref{sec:Dick-alpha}.
We omit the definition of $||f||_\alpha$,
which depends on all partial mixed derivatives
up to the $\alpha$-th order in each variable;
when $s=1$, it is defined by $$\| f\| ^2_\alpha:=
\sum_{i=0}^{\alpha}\Big|\int_0^1 f^{(i)}(x)\, dx\Big|^2+\int_0^1\Big| f^{(\alpha)}(x)\Big| ^2\, dx.$$

\section{Digital net, Discretization and WAFOM}
In \cite{WAFOM}, Saito, Matoba and the first author
introduced Walsh figure of merit (WAFOM) $\wafom(P)$
of a digital net\footnote{See \S\ref{sec:digital} for 
a definition of digital nets; there we use the italic $P$ instead of 
$\calP$ for a digital net, to stress that actually $P$ is a
subspace of a discrete space, while $\calP$ is in a 
continuous space $I^s$.} $P$.
This may be regarded as a simplified special case of 
Dick's $W_\alpha$ with some discretization.
WAFOM satisfies a Koksma-Hlawka type inequality,
and the value $\wafom(P)$ decreases 
in the order $O(N^{-C(\log_2 N)/s+D})$ for some
constant $C, D>0$ independent of $s, N$.
Thus, the order of the convergence is faster than $O(N^{-\alpha})$
for any $\alpha>0$. 

\subsection{Discretization}
Although the following notions are
naturally extended to $\Z/b$ or 
even any finite abelian groups
\cite{SUZUKI2014}, we treat only the case
when base $b=2$ for simplicity. 

Let $\Ftwo:=\{0,1\}=\Z/2$ be the two-element field. 
Take $n$ large enough, and 
approximate the unit interval $I=[0,1)$
by the set of $n$-bit integers 
$I_n:=\Ftwo^n$
through the inclusion 
$I_n \to I$, $x (\mbox{considered as an $n$-bit integer})\mapsto x/2^n + 1/2^{n+1}$.

More precisely, 
we identify the finite set $I_n$
with the set of half open intervals
obtained by partitioning $[0,1)$ into $2^n$ pieces;
namely 
$$
 \calI_n:=\{[i 2^{-n}, (i+1)2^{-n}) \ | \ 0\leq i \leq 2^n-1\}.
$$

\begin{figure}
\caption{$\{0,1\}^3$ is identified with the set of 8 segments $\calI_3$.}
\label{fig:3-digits}
\begin{center}
\includegraphics[width=0.6\textwidth, bb=0 100 800 400]{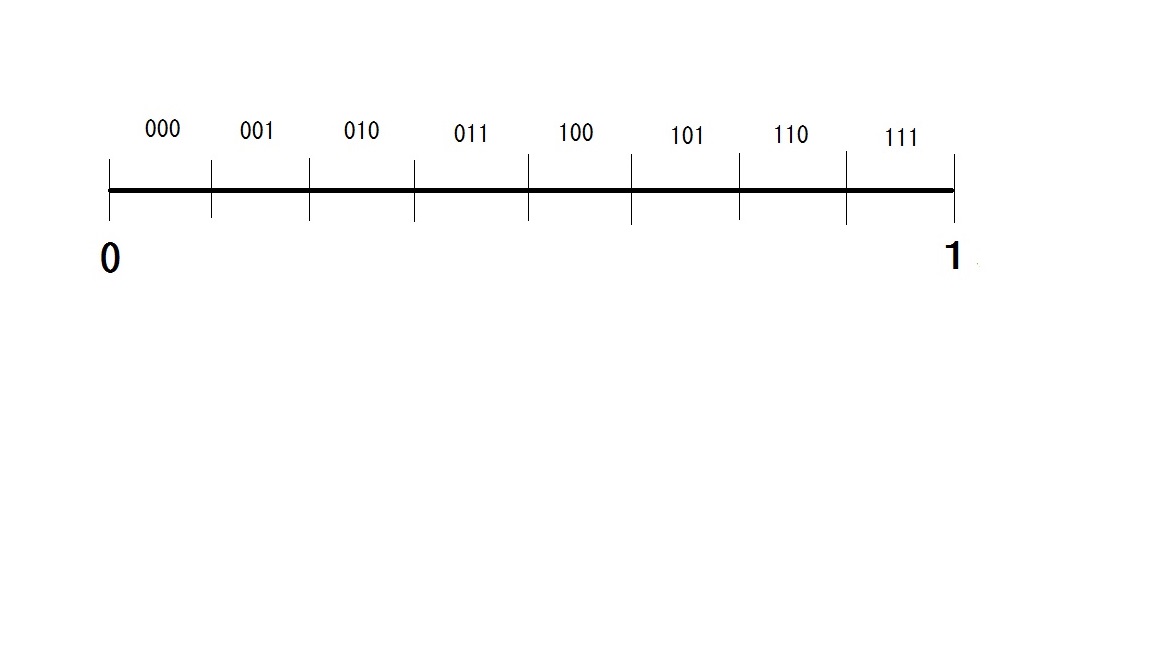}
\end{center}
\end{figure}

\begin{example} In the case $n=3$ and $I_3=\{0,1\}^3$, 
$\calI_3$ is the set of 8 intervals in Figure~\ref{fig:3-digits}.
\end{example}
The $s$-dimensional hypercube $I^s$
is approximated by the set $\calI_n^s$ of $2^{ns}$ hypercubes,
which is identified with $I_n^s=(\Ftwo^n)^s=M_{s,n}(\Ftwo)=:V$.
In sum, 
\begin{definition}
Let $V:=M_{s,n}(\Ftwo)$ be the set of $(s\times n)$-matrices
with coefficients in $\Ftwo=\{0,1\}$. An element $B=(b_{ij})\in V$
is identified with an $s$-dimensional hypercube in
$\calI_{n}^s$, consisting of elements $(x_1, \ldots, x_s)\in \R^s$
where, for each $i$,
the binary expansion of $x_i$ coincides with 
$0.b_{i1}b_{i2}\cdots b_{in}$
up to the $n$-th digit below the decimal point.
By abuse of the language, the notation $B$ is used for 
the corresponding hypercube.
\end{definition}
\begin{example} In the case $n=3$ and $s=2$, for example,
$$
B=
\left(
\begin{array}{c}
100 \\
011
\end{array}
\right)
\mbox{ corresponds to }
 [0.100,0.101)\times 
 [0.011,0.100). 
$$
\end{example}

As an approximation of $f:I^s \to \R$, define
$$f_n:\calI_n^S=V \to \R, \quad B \mapsto f_n(B):=\frac{1}{\Vol(B)}\int_{B} f \rd x$$
by mapping a small hypercube $B$ of edge length $2^{-n}$
to the average of $f$ over this small hypercube.
Thus, $f_n$ is the discretization (with $n$-bit precision) of $f$ by taking the average
over each small hypercube. 

In the following, we do not compute $f_n$, but 
consider as if we are given $f_n$. More precisely saying,
let $\bsx_B$ denote the mid point of the hypercube $B$,
and we approximate $f_n(B)$ by $f(\bsx_B)$. 
For sufficiently large $n$, say, $n=32$, the approximation 
error $|f_n(B)-f(\bsx_B)|$ 
(which we call the discretization error of
$f$ at $B$
)
would be small enough: if $f$ is Lipschitz continuous,
then the error\footnote{
 If $f$ has Lipschitz constant $C$, namely, satisfies $f(x-y)<C|x-y|$, 
then the error is bounded by $C\sqrt{s}2^{-n}$
\cite[Lemma~2.1]{WAFOM}.
} has order $\sqrt{s}2^{-n}$. 

From now on, we assume that $n$ is taken large enough,
so that this discretization error is negligible in 
practice for the QMC integration considered. 
A justification is that we have only finite precision 
computation in digital computers, so a function 
$f$ has discretized domain with some finite precision.
This assumption is somewhat cheating, but seems to work
well in many practical uses.

By definition of the above discretization, we have an equality
$$
\int_{[0,1)^s} f(x)\rd x = 
\frac{1}{|V|}\sum_{B \in V} f_n(B).
$$

\subsection{Discrete Fourier transform}
For $A, B \in V$, we define its inner product by
$$
(A,B):= \mbox{trace} (^tA B) = 
\sum_{1\leq i \leq s, 1\leq j \leq n} a_{ij}b_{ij}\in \Ftwo \quad (\bmod 2).
$$
For a function
$
g: V \to \R, 
$
its \emph{discrete Fourier transform} $\hat{g}:V \to \R$ is
defined by 
$$
\hat{g}(A):= \frac{1}{|V|}
\sum_{B \in V} g(B)(-1)^{(B,A)}.
$$
Thus
$$
\hat{f_n}(0)=\frac{1}{|V|}\sum_{B \in V} f_n(B) =I(f).
$$
\begin{remark} 
The value $\hat{f_n}(A)$ coincides with
the $A$-th \emph{Walsh coefficient} of the function $f$ 
defined as follows.
Let $A=(a_{ij})$. Define an integer $c_i:=\sum_{j=1}^n a_{ij}2^j$
for each $i=1,\ldots,s$.
Then the $A$-th Walsh coefficient of $f$ is defined
as the standard multi-indexed 
Walsh coefficient $\hat{f}_{c_1,\ldots,c_s}$.
\end{remark}
\subsection{Digital nets, and QMC-error in terms of Walsh coefficients}
\label{sec:digital}
\begin{definition}
Let $P\subset V$ be an $\Ftwo$-linear subspace
(namely, $P$ is closed under componentwise addition modulo 2).
Then, $P$ can be regarded as a set of small hypercubes
in $\calI_n^s$, or, a finite point set $\calP \subset I^s$
by taking the mid point of each hypercubes.
Such a point set $\calP$ (or even $P$) is called a digital net with base 2.
\end{definition}
This notion goes back to Sobol' and Niederreiter;
see for example \cite[Definition~4.47]{DICK-PILL-BOOK}.
For such an $\Ftwo$-subspace $P$, let us define 
its perpendicular space\footnote{
 The perpendicular space is called ``the dual space''
in most literatures on QMC and coding theory. However,
in pure algebra, the dual space to a vector space $V$ 
over a field $k$ means $V^*:=\mathrm{Hom}_k(V,k)$,
which is defined without using inner product. In this paper,
we use the term ``perpendicular'' going against
the tradition in this area.
}
by
$$P^\perp:=\{A \in V\ | \ (B,A)=0 \ (\forall B \in P)\}.$$
QMC integration of $f_n$ by $P$ is by definition
\begin{equation}\label{eq:byWalsh}
I(f_n;P):=\frac{1}{|P|}\sum_{B \in P} f_n(B) 
=
\sum_{A \in P^\perp} \hat{f_n}(A),
\end{equation}
where the right equality (called Poisson summation formula)
follows from
$$
\begin{array}{rcl}
\sum_{A \in P^\perp} \hat{f_n}(A)
&=&
\sum_{A \in P^\perp} \frac{1}{|V|}
(\sum_{B \in V} f_n(B)(-1)^{(B,A)}) \\
&=&
\frac{1}{|V|}
\sum_{B \in V} f_n(B)\sum_{A \in P^\perp} (-1)^{(B,A)} \\
&=&
\frac{1}{|V|}
\sum_{B \in P} f_n(B) |P^\perp| \\
&=&
\frac{1}{|P|}
\sum_{B \in P} f_n(B). \\
\end{array}
$$
\subsection{Koksma-Hlawka type inequality by Dick}
From (\ref{eq:byWalsh}), we have a QMC integration error bound
by Walsh coefficients
\begin{equation}\label{eq:boundbyWalsh}
\Error(f_n;P)= |I(f_n;P) - \hat{f_n}(0)| \\
=\left|\sum_{A \in P^\perp-\{0\}} \hat{f_n}(A)\right| 
\leq \sum_{A \in P^\perp-\{0\}} |\hat{f_n}(A)|. \\
\end{equation}
Thus, to bound the error, it suffices to bound $|\hat{f_n}(A)|$.

\begin{theorem}[Decay of Walsh coefficients, \cite{Dick-Decay}]
\label{th:decay}
For an $n$-smooth function $f$, there is a notion of
$n$-norm $||f||_n$ and a constant
$C(s,n)$ independent of $f$ and $A$ with
$$
|\hat{f_n}(A)|\leq C(s,n)||f||_n 2^{-\mu(A)}.
$$
\end{theorem}
(See \cite[Theorem~14.23]{DICK-PILL-BOOK} for a general statement.)
Here, $\mu(A)$ is defined as follows:
\begin{definition}\label{def:Dick-weight}
For $A=(a_{ij})_{1\leq i\leq s, 1\leq j \leq n} \in V$, its Dick weight 
$\mu(A)$ is defined by 
$$
\mu(A):=\sum_{1\leq i \leq s, 1\leq j \leq n} ja_{ij},
$$
where $a_{ij}\in \{0,1\}$ are considered as integers
(without modulo 2).
\end{definition}

\begin{example} In the case of $s=3, n=4$, 
for example, 
$$
A=
\left(
\begin{array}{c}
1001 \\
0111 \\
0010
\end{array}
\right)
\stackrel{ja_{ij}}{\to}
\left(
\begin{array}{c}
1004 \\
0234 \\
0030
\end{array}
\right)
\to
\mu(A)=
\begin{array}{r}
(1+0+0+4)\\
+(0+2+3+4)\\
+(0+0+3+0)
\end{array}
=17.
$$
\end{example}

Walsh figure of merit of $P$ is defined as follows \cite{WAFOM}:
\begin{definition}[WAFOM]\label{def:WAFOM}
Let $P\subset V$. WAFOM of $P$ is defined by
$$
\wafom(P):=\sum_{A \in P^\perp-\{0\}}2^{-\mu(A)}.
$$
\end{definition}
By plugging this definition and Dick's Theorem~\ref{th:decay} 
into (\ref{eq:boundbyWalsh}), we have an inequality of 
Koksma-Hlawka type:
\begin{equation}\label{eq:boundbyWAFOM}
\Error(f_n;P)
\leq C(s,n)||f||_n \wafom(P).
\end{equation}


\subsection{A toy experiment on $\wafom(P)$}
We shall see how WAFOM works
for a toy case of $n=3$-digit precision and $s=1$ dimension.
In Figure~\ref{fig:3-digits}, the unit interval $I$
is divided into 8 intervals, each of which corresponds
to a $(1\times 3)$-matrix in $\Ftwo^3=V$.
Table\ref{tab:toy} lists the seven subspaces of dimension 2,
selection of four of them, and
their WAFOM and QMC error for the integrand $f(x)=x,x^2$ and $x^3$.
\begin{table}
\caption{Toy examples for WAFOM for 3-digit discretization for 
integrated $x, x^2$ and $x^3$}\label{tab:toy}
$$
\begin{array}{rccccccl}
V
=\{
000 & 001 & 010 & 011 & 100 & 101 & 110 & 111 
\} \\
(100)^\perp
=\{
000 & 001 & 010 & 011 & \phantom{000} & \phantom{000} & \phantom{000} & \phantom{000} 
\} \\
(010)^\perp
=\{
000 & 001 & \phantom{000}& \phantom{000} & 100 & 101 & \phantom{000} & \phantom{000} 
\} \\
(110)^\perp
=\{
000 & 001 & \phantom{000}& \phantom{000} & \phantom{000} & \phantom{000} & 110 & 111 
\} \\
(001)^\perp
=\{
000 & \phantom{000} & 010 & \phantom{000} & 100 & \phantom{000} & 110 & \phantom{000}
\} \\
(101)^\perp
=\{
000 & \phantom{000} & 010 & \phantom{000} & \phantom{000} & 101 & \phantom{000} & 111 
\} \\
(011)^\perp
=\{
000 & \phantom{000} & \phantom{000} & 011 & 100 & \phantom{000} & \phantom{000} & 111 
\} \\
(111)^\perp
=\{
000 & \phantom{000} & \phantom{000} & 011 & \phantom{000}& 101 & 110 & \phantom{000} 
\} \\
\end{array}
$$
\begin{center}
\begin{tabular}{|c|c|c||c|c|c|}
$P$
& 
$\begin{array}{l}\mu(A) \mbox{ for}\\ A\in P^\perp
\setminus 0
\end{array}$
&WF($P$)
& Error for $x$ &Error for $x^2$ &Error for $x^3$ \\
\hline
$V$ & $\emptyset$ &0 & 0 & $-$0.0013 & $-$0.0020 \\
$001^\perp$ & 0+0+3 &$2^{-3}$ &$-$0.0625 & $-$0.0638 & $-$0.0637 \\
$101^\perp$ & 1+0+3 &$2^{-4}$& 0 & $-$0.0299 & $-$0.0449 \\
$011^\perp$ & 0+2+3 &$2^{-5}$& 0 & +0.0143 & +0.0215 \\
$111^\perp$ & 1+2+3&$2^{-6}$ & 0 & $-$0.0013 & $-$0.0137 \\
\end{tabular}
\end{center}
\end{table}
The first line in Table~\ref{tab:toy} shows the 8-element
set $V=\Ftwo^3$, corresponding to the 8 intervals in 
Figure~\ref{fig:3-digits}.
The next line $(100)^\perp$ denotes the 2-dimensional 
subspace of $V$ consisting of the elements perpendicular
to $(100)$, that is, the four vectors whose first digit is 0.
In the same manner, all 2-dimensional subspaces of $V$
are listed. The last one is $(111)^\perp$, 
consisting of the four vectors $(x_1,x_2,x_3)$ with $x_1+x_2+x_3=0 (\bmod 2)$.

Our aim is to decide which is the best (or most ``uniform'') 
among the seven 
2-dimensional sub-vector spaces for QMC integration.
Intuitively, $(100)^\perp$ is not a good choice since 
all the four intervals cluster in $[0,1/2]$.
Similarly, we exclude $(010)^\perp$ and $(110)^\perp$.
We compare the remaining four candidates by two methods:
computing WAFOM, and computing QMC
integration errors with test integrand functions $x, x^2$ and $x^3$.

The results are shown in the latter part of Table~\ref{tab:toy}.
The first line corresponds to the case of $P=V$.
Since $P^\perp -\{0\}$ is empty, $\wafom(P)=0$.
For the remaining four cases $P=(x_1, x_2, x_3)^\perp$, note that
$\{(x_1, x_2, x_3)^\perp\}^\perp=\{(000), (x_1,x_2,x_3)\}$
and $P^\perp -\{0\}=\{(x_1,x_2,x_3)\}$, thus we have
$
\wafom(P)=2^{-\mu((x_1,x_2,x_3))}.
$
The third column in the latter table 
shows WAFOM for five different choices of $P$. 
The three columns ``Error for $x^i$'' with $i=1,2,3$
show the QMC integration error by $P$ for integrating $x^i$
over $[0,1]$. We used the mid point of each segment (of length 1/8)
to evaluate $f$. Thus, the listed errors include both
the discretization errors and QMC-integration errors for $f_n$.
For the first line, $P=V$ implies no QMC integration error for $f_n$
($n=3$), so the values show the discretization error exactly. 
The error bound (\ref{eq:boundbyWAFOM}) 
is proportional to $\wafom(P)$ for a fixed integrand. The table shows 
that, for these test functions, the actual errors
are well reflected in WAFOM values.

Here is a loose interpretation of $\wafom(P)$.
For an $\Ftwo$-linear $P$, 
\begin{itemize}
\item $A \in P^\perp\setminus \{0\}$ is a linear relation satisfied by $P$.
\item $\mu(A)$ measures ``complexity'' of $A$.
\item $\wafom(P)=\sum_{A \in P^\perp \setminus \{0\}}2^{-\mu(A)}$
is small if all relations have high complexity, 
and hence $P$ is close to ``uniform.''
\end{itemize}
The weight $j$ in the sum $\sum ja_{ij}$ in the definition of $\mu(A)$
denotes that the $j$-th digit below the decimal point is counted
with complexity $2^{-j}$.

\section{Point sets with low WAFOM values}

\subsection{Existence and non-existence of low WAFOM point sets}
\begin{theorem} \label{th:existence}
There are absolute (i.e.\ independent of $s, n$ and $d$) positive constants $C, D, E$ such that
for any positive integer $s, n$ and $d\geq 9s$,
there exists a $P \subset V$ of $\Ftwo$-dimension $d$
(hence cardinality $N=2^d$) satisfying
$$
\begin{array}{rlc}
\wafom(P)&\leq &E\cdot 2^{-Cd^2/s + Dd} 
= E\cdot N^{-C\log_2 N/s + D}.
\end{array}
$$
\end{theorem}
Since the exponent $-C\log_2 N/s + D$ goes to $-\infty$
when $N \to \infty$, this shows that there exist
point sets with ``higher order convergence''
having this order of WAFOM.
There are two independent proofs:
M-Yoshiki \cite{YOSHIKI-MATSUMOTO}
shows the positivity of the probability to have low-WAFOM
point sets under a random choice of its basis (hence non-constructive),
and K.Suzuki \cite{SUZUKI2012} shows a construction 
using Dick's interleaving method \cite[\S15]{DICK-PILL-BOOK}
for Niederreiter-Xing sequence \cite{niederreiter:xing}.
Suzuki \cite{SUZUKI2014} generalizes 
\cite{YOSHIKI-MATSUMOTO} and \cite{YOSHIKI2014} for arbitrary base $b$.
Theorem~\ref{th:existence} is similar to the Dick's construction
of point sets with $W_\alpha(\calP)=O(N^{-\alpha}(log N)^{s\alpha})$
for arbitrary high $\alpha \geq 1$, but there seems no 
implication between his result and this theorem.

On the other side, Yoshiki \cite{YOSHIKI2014} proved 
the following theorem that
the order of the exponent $d^2/s$ is sharp, namely, WAFOM
can not be so small:
\begin{theorem} 
$ $ 
Let $C'>1/2$ be any constant. 
For any positive integer $s, n$ and $d\geq s\times(\sqrt{C'+1/16}+3/4)/(C'-1/2)$,
any linear subspace $P \subset V$ of $\Ftwo$-dimension $d$
satisfies
$$
\begin{array}{rlc}
\wafom(P)&\geq &2^{-C'd^2/s}. \\
\end{array}
$$
\end{theorem}

\subsection{An efficient computation method of WAFOM}
Since $P$ is intended for a QMC integration
where the enumeration of $P$ is necessary, $|P|=2^{\dim_{\Ftwo}P}$
can not be huge.
On the other hand, $|V|=2^{ns}$ would be huge, say, for $n=32$ and $s>2$.
Since $\dim_{\Ftwo} P + \dim_{\Ftwo} P^\perp = \dim_{\Ftwo} V$,
$|P^\perp|$ must be huge. Thus, a direct computation of 
$\wafom(P)$ using Definition~\ref{def:WAFOM} would be too costly.
In \cite{WAFOM}, the following formula is given by a Fourier inversion.
Put $B=(b_{i,j})$, then we have
$$
\wafom(P)=
\frac{1}{|P|}
\sum_{B \in P}
\left\{ 
\prod_{1\leq i \leq s, 1 \leq j\leq n} [(1+(-1)^{b_{i,j}}2^{-j})] -1 
\right\}.
$$
This is computable in $O(nsN)$ steps of arithmetic
operations in real numbers, where $N=|P|$.
Compared with most of other discrepancies, this is 
relatively easily computable. This allows us to do a
random search for low-WAFOM point sets.

\begin{remark}
\begin{enumerate}
\item The above equality holds only for an $\Ftwo$-linear $P$. 
Since the left hand side is non-negative, so is the
right sum in this case. 
It seems impossible to define WAFOM for
a general point set by using this formula, since for a general (i.e.\ non-linear) $P$, 
the sum at the right hand side is sometimes negative and thus will never give
a bound on the integration error.
\item
The right sum may be interpreted as the QMC integration of 
a function (whose definition is given in the right hand side
of the equality) by $P$. The integration of the function 
over total space $V$ is zero. Hence, the above equality indicates that,
to have a best $\Ftwo$-linear $P$ from the viewpoint of WAFOM,
it suffices to have a best $P$ for QMC integration 
for a single specified function. This is in contrast to
the definition of star-discrepancy, where all the rectangle
characteristic functions are used as 
the test functions, and the supremum 
of their QMC integration errors is taken.
\item Harase-Ohori\cite{HARASE-OHORI} gives a method 
to accelerate this computation by a factor of 30,
using a look-up table. Ohori-Yoshiki\cite{OHORI-YOSHIKI}
gives a faster and simpler method to compute a good 
approximation of WAFOM, using that Walsh coefficients
of exponential function approximates the Dick weight $\mu$.
More precisely, $\wafom(P)$ is well-approximated by 
the QMC-error of the function $\exp(-2\sum_{i=1}^s x_i)$, 
whose value is easy to evaluate in modern CPUs.
\end{enumerate}
\end{remark}

\section{Experimental results}
\subsection{Random search for low WAFOM point sets}
We fix the precision 
$n=30$. We consider two cases of the dimension $s=4$ and $s=8$.
For each $d=8,9,10,\ldots, 16$, we generate 
$d$-dimensional subspace $P \subset V=(\Ftwo^{30})^s$
10000 times, by the uniformly random choice of 
$d$ elements as its basis.
Let $P_{d,s}$ be the point set with the lowest WAFOM among them.
For the comparison, $Q_{d,s}$ be the 
point set of the 100th lowest WAFOM.

\subsection{Comparison of QMC rules by WAFOM}
\begin{figure}
 \caption{WAFOM values for: (1) best WAFOM among 10000, (2)
 the 100th best WAFOM, (3) Niederreiter-Xing, (4)
 Sobol',
 of size $2^d$ with $d=8,9,\ldots,16$.
 The vertical axis is for $\log_2$ of their WAFOM, and the horizontal
 for $\log_2$ of the size of point sets. The left figure is for dimension $s=4$,
 the right $s=8$.
 }\label{fig:WAFOM}
 \begin{center}
  \begin{tabular}{c}

   \begin{minipage}{0.5\hsize}
    \begin{center}
    \includegraphics[width=1.0\textwidth]{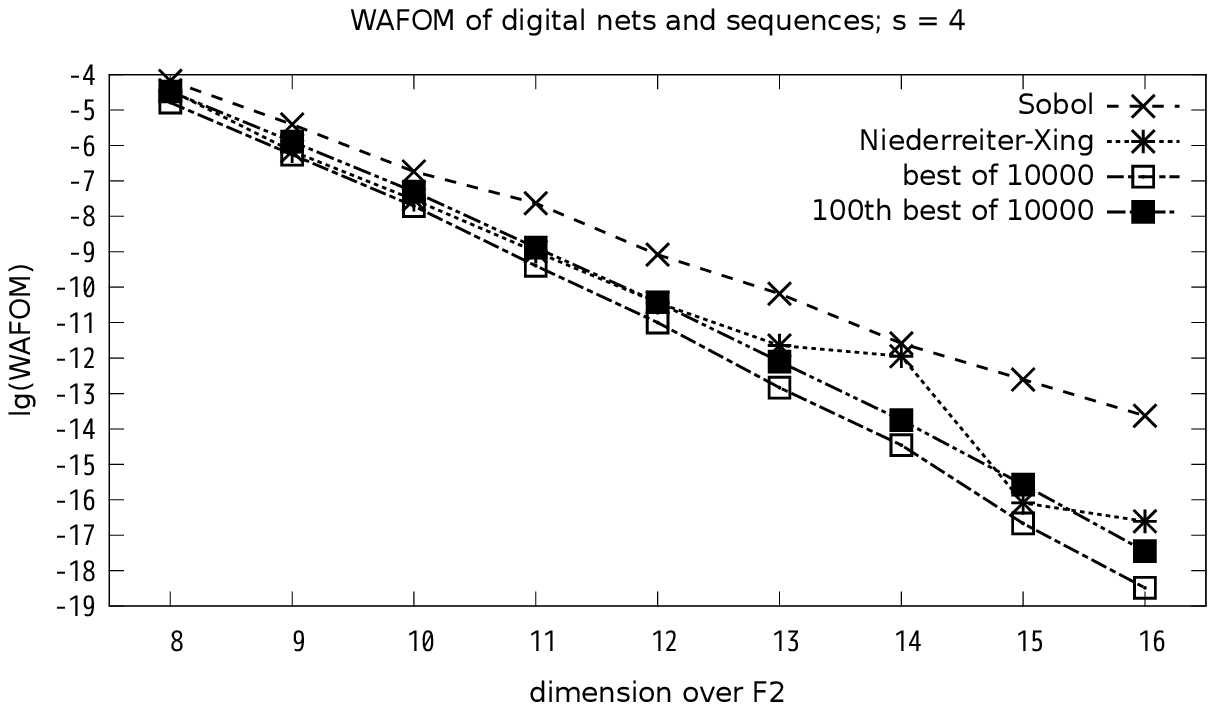}
    \end{center}
   \end{minipage}

   \begin{minipage}{0.5\hsize}
    \begin{center}
    \includegraphics[width=1.0\textwidth]{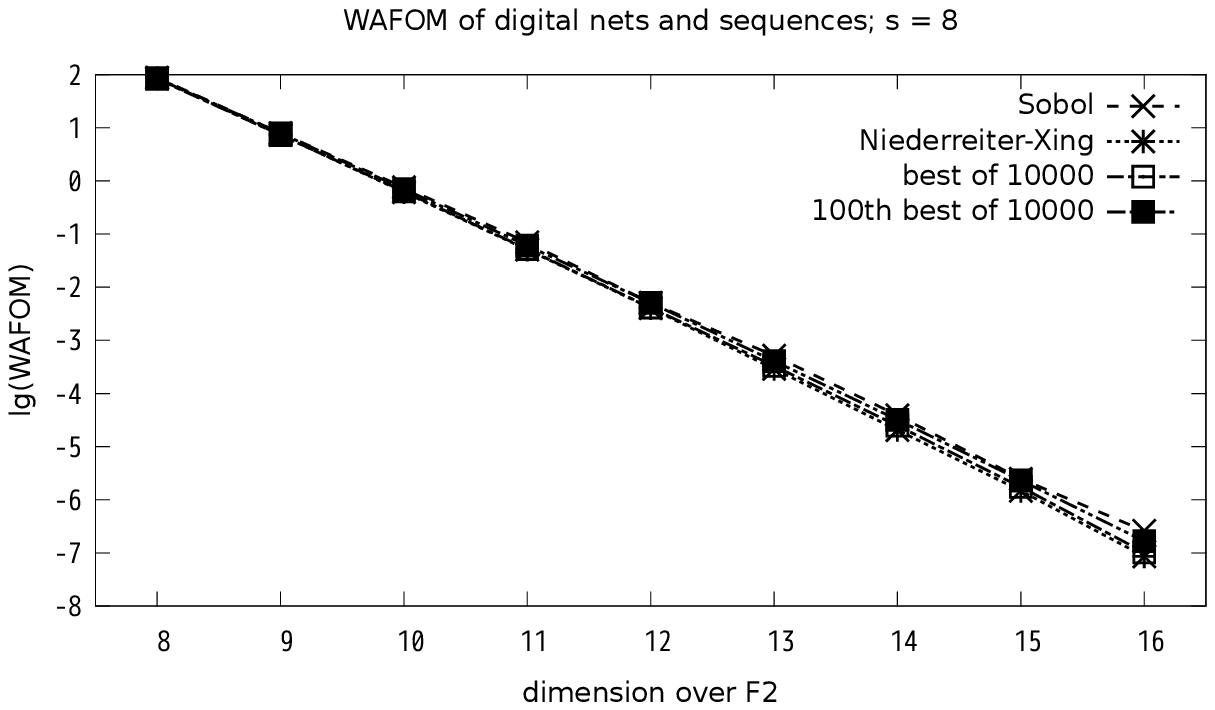}
    \end{center}
   \end{minipage}
  \end{tabular}
 \end{center}
\end{figure}

For a comparison, we use two other QMC quadrature rules, namely, Sobol' sequence 
improved by Joe and Kuo \cite{JOE-KUO},
and Niederreiter-Xing sequence (NX) implemented by Pirsic \cite{PIRSIC}
and by Dirk Nuyens \cite[item {\tt nxmats}]{NUYENS} (downloaded from the latter).
Figure~\ref{fig:WAFOM} shows the WAFOM values for
these four kinds of point sets, with size $2^8$ to $2^{16}$.
For $s=4$, Sobol' has largest WAFOM value,
while NX has small WAFOM comparable to the 100th best 
$Q_{d,s}$ selected by WAFOM. In $d=14$, NX has much larger WAFOM
than that of $Q_{14,s}$, while in $d=15$ the converse occurs.
Note that this seems to be reflected in the following experiments.
For $s=8$, the four kinds of point sets show small differences 
in values of their WAFOM.
Indeed, NX has smaller WAFOM value than the best point set
among randomly generated 10000 for each $d$, while Sobol'
has larger WAFOM values. A mathematical analysis on this
good grade of NX would be interesting.

\subsection{Comparison by numerical integration}
In addition to the above four kinds of QMC rules,
Monte Carlo method is used for comparison (using 
Mersenne Twister \cite{MT} pseudorandom number generator).
For the test functions, we use 6 Genz functions \cite{GENZ1987}:
\begin{description}
\item[Oscillatory]
$f_1(\bsx)=\cos(2\pi u_1+\sum_{i=1}^s a_i x_i)$,
\item[Product Peak]
$f_2(\bsx)=\prod_{i=1}^s[1/(a_i^2+(x_i-u_i)^2)]$,
\item[Corner Peak]
$f_3(\bsx)=(1+\sum_{i=1}^s a_i x_i)^{-(s+1)}$
\item[Gaussian]
$f_4(\bsx)=\exp(-\sum_{i=1}^s a_i^2(x_i-u_i)^2)$
\item[Continuous]
$f_5(\bsx)=\exp(-\sum_{i=1}^s a_i|x_i-u_i|)$
\item[Discontinuous]
$
f_6(\bsx)=
\begin{cases}
0 & \mbox{ if } x_1>u_1 \mbox{ or } x_2>u_2, \\
\exp(\sum_{i=1}^s a_i x_i)) & \mbox{otherwise.}
\end{cases}
$
\end{description}
This selection is copied from \cite[P.91]{NOVAK-RITTER} \cite{HARASE-OHORI}. The parameters
$a_1,\ldots,a_s$ are selected so that (1) they are in an arithmetic
progression (2) $a_s=2a_1$ (3) the average of $a_1,\ldots,a_s$
coincides with the average of $c_1,\ldots,c_{10}$ in
\cite[Equation (10)]{NOVAK-RITTER}
for each test function.
The parameters $u_i$ are generated randomly by \cite{MT}.

\begin{figure}
 \caption{QMC integration errors for (1) best WAFOM among 10000, (2)
 the 100th best WAFOM, (3) Niederreiter-Xing, (4) Sobol', (5) Monte
 Carlo, using six Genz functions on the 4-dimensional unit cube.
 The vertical axis is for $\log_2$ of the errors, and the horizontal
 for $\log_2$ of the size of point sets. The error
 is the mean square error for 100 randomly digital shifted point sets.
 }\label{fig:4dim}
 \begin{center}
  \begin{tabular}{c}

   \begin{minipage}{0.5\hsize}
    \begin{center}
    \includegraphics[width=1.0\textwidth]{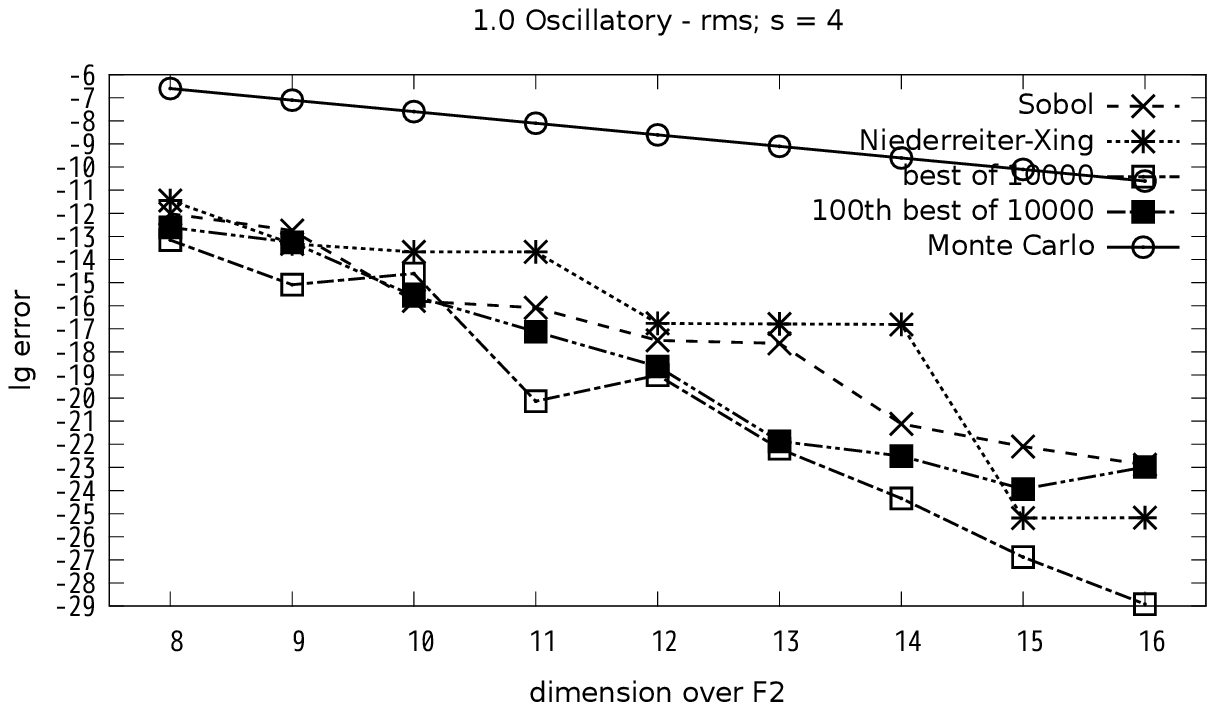}
    \end{center}
   \end{minipage}

   \begin{minipage}{0.5\hsize}
    \begin{center}
    \includegraphics[width=1.0\textwidth]{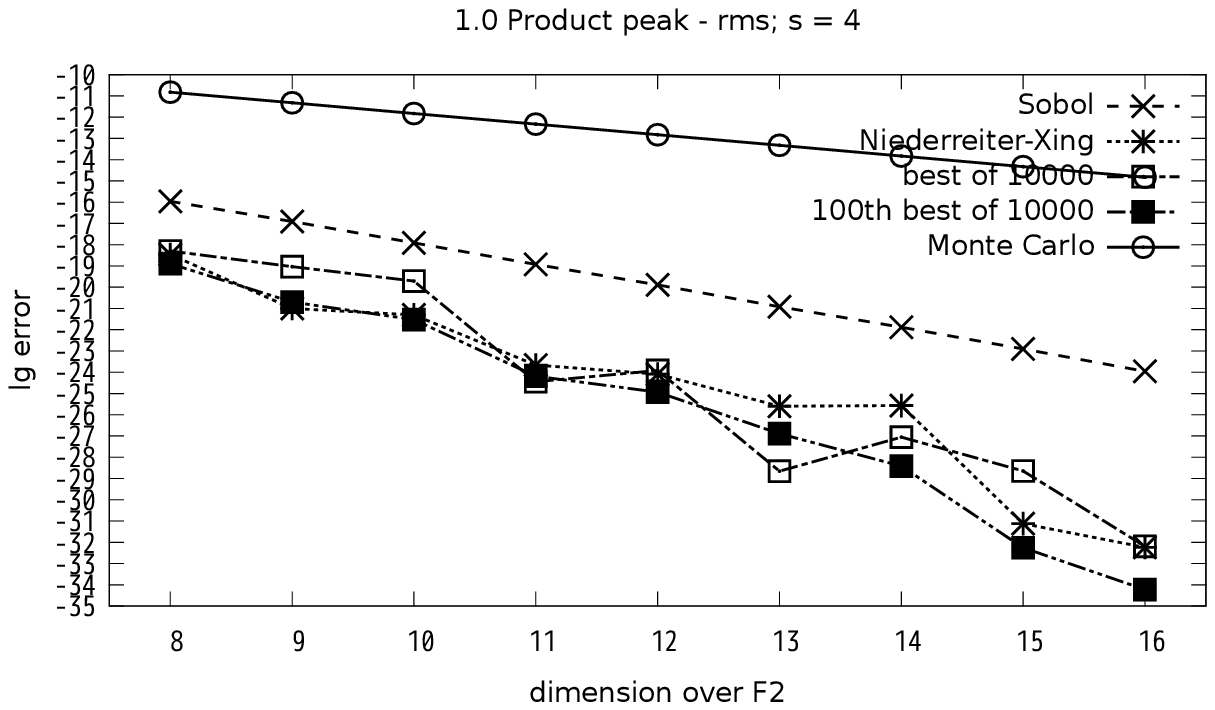}
    \end{center}
   \end{minipage}
   \\
   \begin{minipage}{0.5\hsize}
    \begin{center}
    \includegraphics[width=1.0\textwidth]{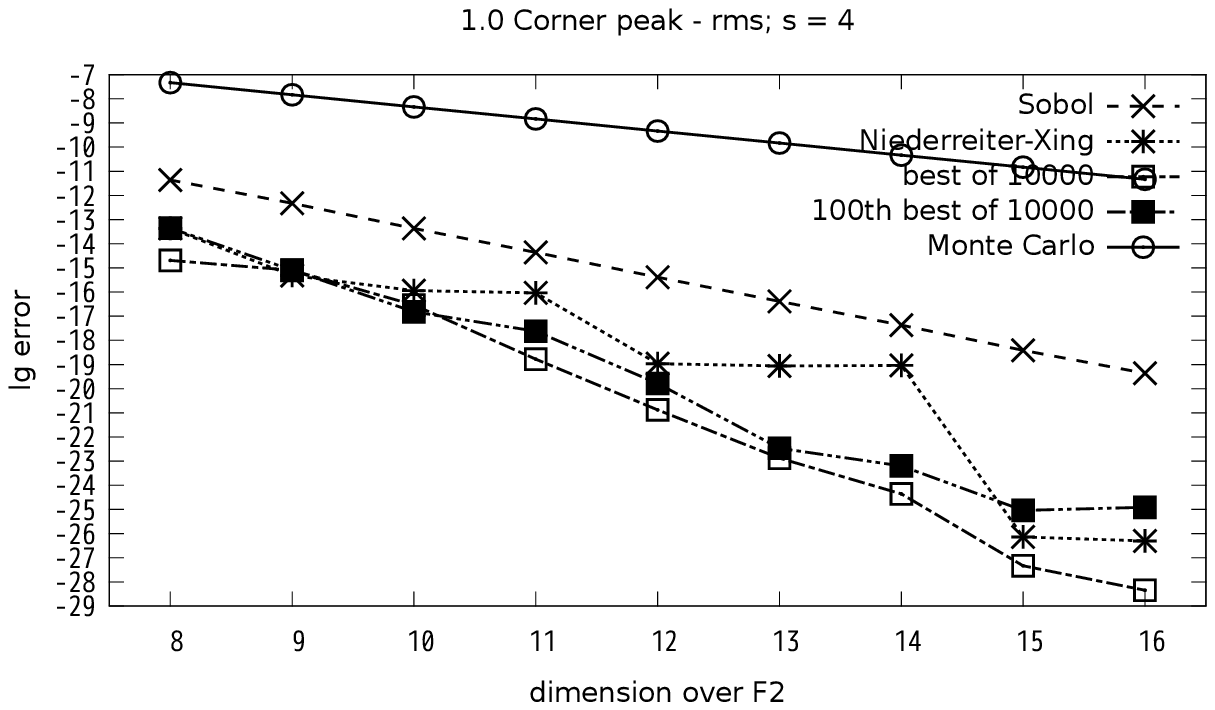}
    \end{center}
   \end{minipage}

   \begin{minipage}{0.5\hsize}
    \begin{center}
    \includegraphics[width=1.0\textwidth]{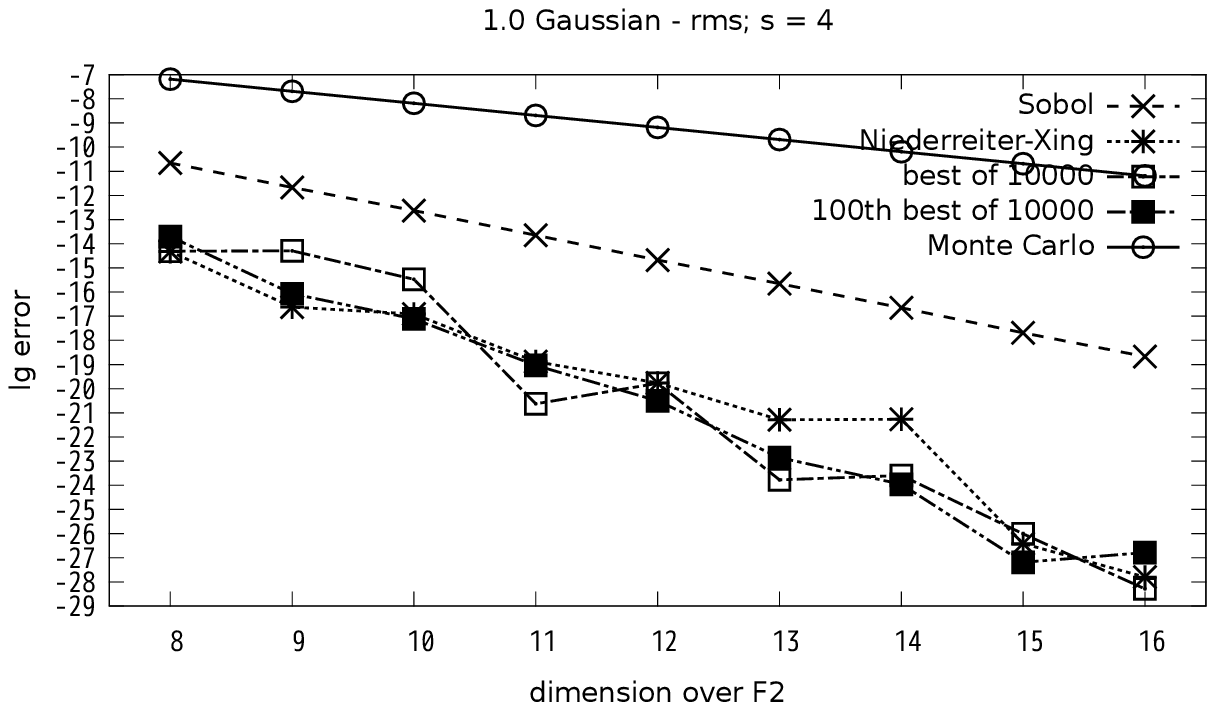}
    \end{center}
   \end{minipage}
   \\
   \begin{minipage}{0.5\hsize}
    \begin{center}
    \includegraphics[width=1.0\textwidth]{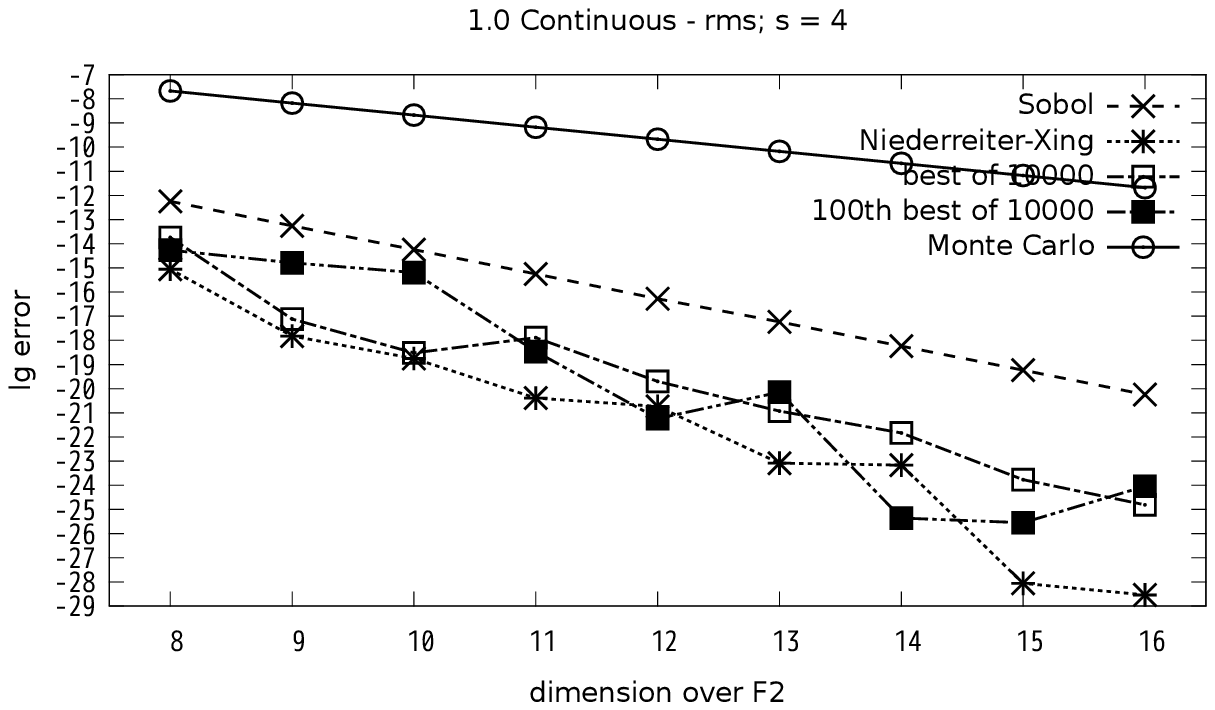}
    \end{center}
   \end{minipage}

   \begin{minipage}{0.5\hsize}
    \begin{center}
    \includegraphics[width=1.0\textwidth]{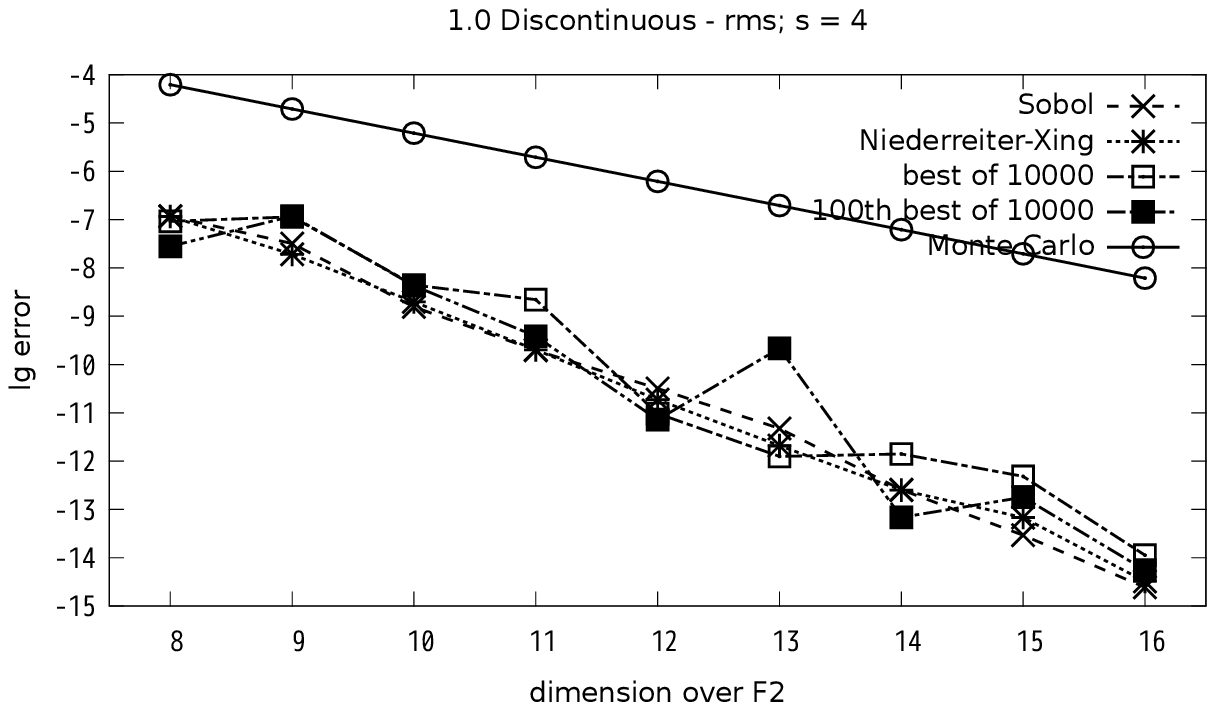}
    \end{center}
   \end{minipage}
   \\
  \end{tabular}
 \end{center}
\end{figure}
Figure~\ref{fig:4dim} shows the QMC integration errors
for six test functions with five methods, for dimension $s=4$. 
The error for Monte Carlo is
of order $N^{-1/2}$. The best WAFOM point sets (WAFOM) and 
Niederreiter-Xing (NX) are comparable.
For the function Oscillatory, where its higher derivatives
grow relatively slowly, WAFOM point sets perform better
than NX and Sobol', and the convergence rate seems
of order $N^{-2}$. For Product peak and Gaussian, WAFOM and NX are
comparable; this coincides with the fact that higher derivatives
of these test functions rapidly grow, but still we observe
convergence rate $N^{-1.6}$.
For Corner peak, WAFOM performs better than NX.
It is somewhat surprising that
the convergence rate is almost $N^{-1.8}$ for WAFOM point sets.
For Continuous, NX performs better than WAFOM.
Since the test functions are not differentiable, 
$||f||_n$ is unbounded 
and hence the inequality (\ref{eq:boundbyWAFOM})
has no meaning. Still, for Continuous, the convergence
rate of WAFOM is almost $N^{-1.2}$. For Discontinuous,
NX and Sobol' perform better than WAFOM.
Note that except Discontinuous, the large/small value of 
WAFOM of NX for $d=14, 15$ observed in the left of 
Figure~\ref{fig:WAFOM} seems to be reflected in
the five graphs.

We conducted similar experiments for $s=8$ dimension,
but we omit the results, since their difference in WAFOM
is small, and the QMC rules show not much difference.
We report that
still we observe
convergence rate with $N^{-\alpha}$ with $\alpha>1.05$
for the five test functions except Discontinuous, for WAFOM selected points
and NX. 
\begin{remark}
\begin{enumerate}
\item
Convergence rate for the integration error is even faster
than that of WAFOM values, for WAFOM selected point sets
and NX for $s=4$, while Sobol' sequence converging
with rate $N^{-1}$. We feel that these go against
our intuition, so checked the code and compared with MC.
We do not know why NX and WAFOM work so well.
\item 
As a referee pointed out, it is hard to observe 
converging rate $N^{-C\log_2 N/s + D}$ in 
Theorem~\ref{th:existence} from the graphs.
\end{enumerate}
\end{remark}

\section{WAFOM versus other figure of merits}
Niederreiter's $t$-value \cite{niederreiter:book}
is a most established figure of merit of a digital net.
Using test functions, we compare the effect of $t$-value
and WAFOM for QMC integration.

\subsection{$t$-value}
Let $\calP\subset I^S=[0,1)^s$ be a finite set of
cardinality $2^m$.
Let $n_1, n_2, \ldots, n_s \geq 0$ be integers.
Recall that $\calI_{n_i}$ is the set of
$2^{n_i}$ intervals partitioning $I$.
Then, $\prod_{i=1}^s\calI_{n_i}$ is 
a set of $2^{n_1+n_2+\cdots + n_s}$
intervals.
We want to make the QMC integration error 0
in computing the volume of every such interval.
A trivial bound is 
$n_1+n_2+\cdots + n_s\leq m$, since at least
one point must fall in each interval.
The point set $\calP$ is called a $(t,m,s)$-net
if the QMC integration error for each interval 
is zero, for any tuple $(n_1,\ldots,n_s)$
with 
$$n_1+n_2+\cdots + n_s\leq m-t.$$
Thus, smaller $t$-value is more preferable.

\subsection{Experiments on WAFOM versus $t$-value}
We fix the dimension $s=4$ and the precision $n=32$,
and generate $10^6$ ($\Ftwo$-linear) point sets of cardinality 
$2^{12}$ by uniform random choices of their $\Ftwo$ basis consisting of 
12 vectors.
We sort these $10^6$ point sets, according to 
their $t$-values. It turns out that 
$3 \leq t \leq 12$, and the frequency of the point sets
for a given $t$-value
is as follows.
\begin{center}
\begin{tabular}{|c||c|c|c|c|c|c|c|c|c|c|}
\hline
$t$ & 3 & 4 & 5 & 6 & 7 & 8 & 9 & 10& 11& 12 \\
freq.
& 63 & 6589 & 29594 & 32403 & 18632 & 
8203 & 2994 & 1059 & 365 & 98 \\
\hline
\end{tabular}
\end{center}
Then, we sort the same $10^6$ point sets by WAFOM.
We categorize them into 10 classes from the
smallest WAFOM, so that $i$-th class has
the same frequency with the $i$-th class
by $t$-value. Thus, the same $10^6$ point sets
are categorized in two ways.
For a given test integrand function,
compute the mean square error of QMC integral
in each category, for those graded by $t$-value
and those graded by WAFOM.

\begin{figure}
 \caption{Left: Hellekalek's function 
     $ f(\bsx)=
     (x_1^{1.1}-\frac{1}{1+1.1})
     (x_2^{1.7}-\frac{1}{1+1.7})
     (x_3^{2.3}-\frac{1}{1+2.3})
     (x_4^{2.9}-\frac{1}{1+2.9})
     $, 
     right: Hamukazu's function 
      $f(\bsx)=
      2^4
      \{5x_1\}
      \{7x_2\}
      \{11x_3\}
      \{13x_4\}
      $,
      where $\{x\}:=x-[x]$.   
  Horizontal axis for category, vertical for the $\log_2$ of
 error. $\Box$:WAFOM, $+ \hspace{-0.75 em}\times$:$t$-value.
 }\label{fig:versus}
 \begin{center}
  \begin{tabular}{c}

   \begin{minipage}{0.5\hsize}
    \begin{center}
    \includegraphics[width=1.0\textwidth]{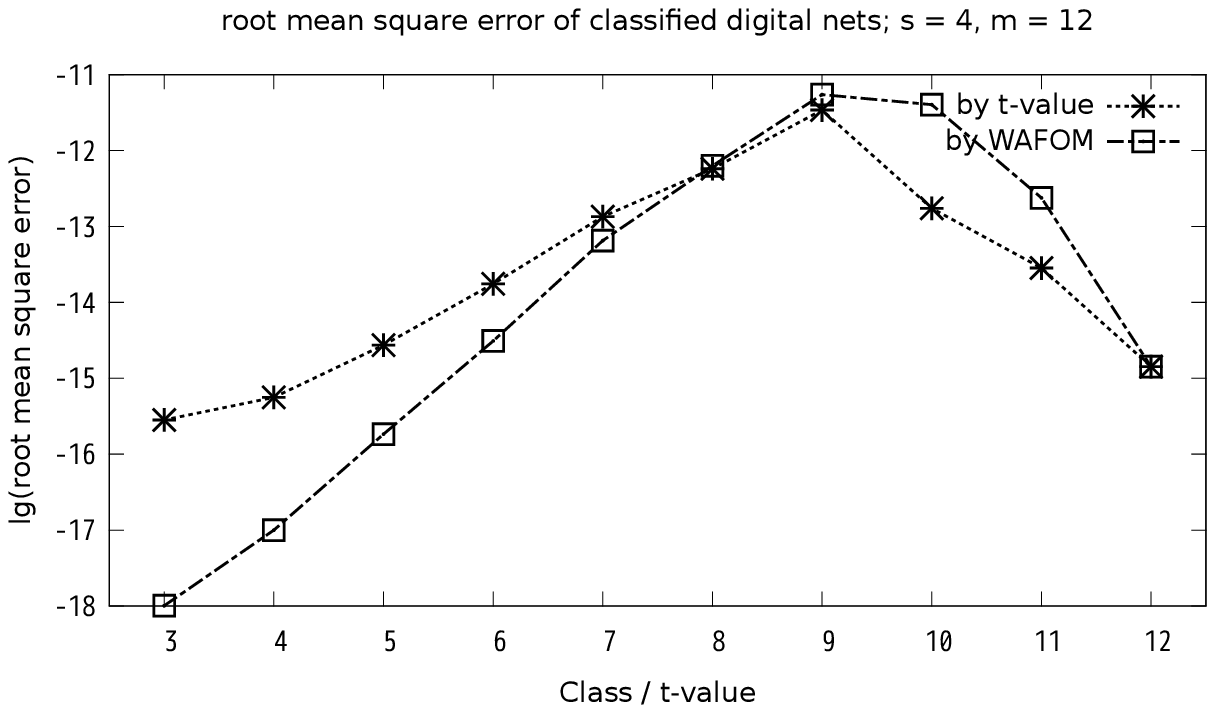}
    \end{center}
   \end{minipage}

   \begin{minipage}{0.5\hsize}
    \begin{center}
    \includegraphics[width=1.0\textwidth]{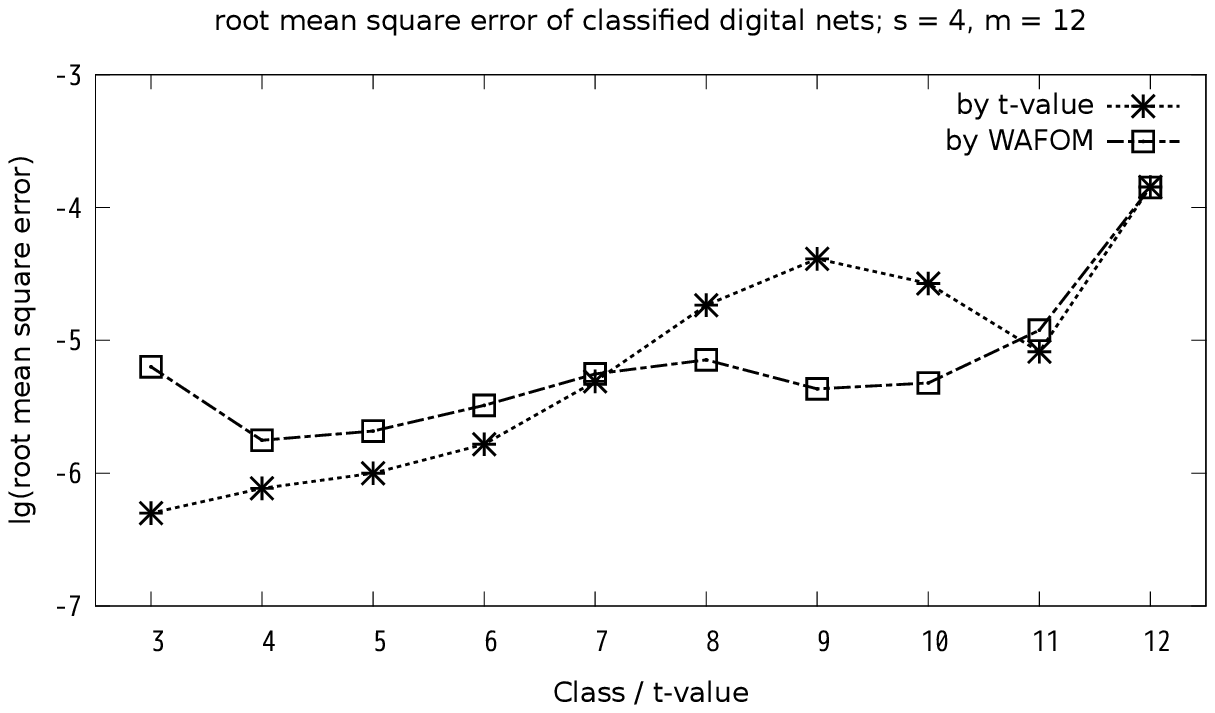}
    \end{center}
   \end{minipage}
  \end{tabular}
 \end{center}
\end{figure}
Figure~\ref{fig:versus} shows $\log_2$ of the mean square integration
error, for each category corresponding to $3\leq t \leq 12$ for
$t$-value ($+ \hspace{-0.75 em}\times$), and for the category
sorted by WAFOM value ($\Box$).
The smooth test function in the left hand side 
comes from Hellekalek \cite{HELEKALLEK1998},
and the non-continuous function in the
right hand side was communicated from Kimikazu Kato
(refered to as ``Hamukazu'' according to his established twitter handle).
From the left figure, for $t=3$,
the average error for the best 63 point sets with 
the smallest $t$-value 3 is much larger 
than the average from the best 63 point sets selected by WAFOM.
Thus, the experiments show that for this test function,
WAFOM seems to work better than $t$-value in selecting
good point set. We have no explanation 
why the error decreases for $t\geq 9$.
In the right figure, for
Hamukazu's non-continuous test function, $t$-value
works better in selecting good points.

Thus, it is expected that digital nets that have small $t$-value 
and small WAFOM would work well for smooth
functions and robust to non-smooth functions.
Harase \cite{HARASE-HYBRID} noticed that
Owen linear scrambling \cite[\S13]{DICK-PILL-BOOK}\cite{OWEN-SCRAMBLING}
preserves $t$-value, 
but changes WAFOM.
Starting from a Niederreiter-Xing sequence with small $t$, 
he applied Owen linear scrambling 
to find a point set 
with low WAFOM and small $t$-value. 
He obtained good results for wide range of integrands.

\subsection{Dick's $\mu_\alpha$, and non-discretized case}\label{sec:Dick-alpha}
Let $\alpha>0$ be an integer. 
For $A \in M_{S,n}(\Ftwo)$, 
the Dick's $\alpha$-weight $\mu_\alpha(A)$
is defined as follows. 
It is a part of summation appeared in Definition~\ref{def:Dick-weight} of $\mu(A)$:
the sum is taken up to $\alpha$ nonzero entries from the right
in each row.

\begin{example} Suppose $\alpha=2$.
$$
A=
\begin{array}{c}
1001 \\
0111 \\
0010
\end{array}
\stackrel{ja_{ij}}{\to}
\begin{array}{c}
\mathbf{1}00\mathbf{4} \\
02\mathbf{3}\mathbf{4} \\
00\mathbf{3}0
\end{array}
\to
\mu_\alpha(A)=
\begin{array}{r}
(1+0+0+4)\\
+(0+0+3+4)\\
+(0+0+3+0)
\end{array}
=15.
$$
\end{example}
For $\Ftwo$-linear $P\subset M_{S,n}(\Ftwo)$,
\begin{equation}\label{eq:W-alpha}
W_\alpha(P):=\sum_{A \in P^\perp-\{0\}}2^{-\mu_\alpha(A)}.
\end{equation}
To be precise, we need to take $n \to \infty$, as follows.
We identify $I=[0,1]$ with 
the product $W:=\Ftwo^\N$ via binary fractional expansion
(neglecting a measure-zero set). 
Let $K:=\Ftwo^{\oplus \N}\subset W$
be the subspace consisting of vectors with finite number of nonzero components
(this is usually identified with $\N \cup \{0\}$ via binary expansion
and reversing the digits).
We define inner product $W \times K \to \Ftwo$ as usual.
Then, for a finite subgroup $P \subset W^s$, its perpendicular
space $P^\perp \subset K^s$ is defined and is countable.
For $A\in K^s$, $\mu_\alpha(A)$ is analogously defined, 
and the right hand side of (\ref{eq:W-alpha}) is 
absolutely converging. Dick \cite{Dick-Decay} proved
$$
\Error(f;P) \leq 
C(s,\alpha)||f||_\alpha W_\alpha(P),
$$
and constructed a sequence of $P$
with $W_\alpha(P)=O(N^{-\alpha}(\log N)^{S\alpha})$
called higher order digital nets.
(See \cite{DICK-PILL-BOOK} for a comprehensive explanation.)
Existence results
and search algorithms for higher order
polynomial lattice rules are studied in 
\cite{BARDEAUX-DICK} \cite{DICK-KRITZER-PILL-SCH}.

WAFOM is an $n$-digit
discretized version of $W_\alpha$ where $\alpha=n$.
WAFOM loses freedom to choose $\alpha$, 
but it might be a merit since we do not need to choose $\alpha$.
\begin{remark} In Dick's theory, $\alpha$ is fixed. 
In fact, setting $\alpha=\log N$ does not
yield useful bound, since
$C(s,\log N)W_{\log N}(P) \to \infty \ (N \to \infty)$.
\end{remark}
%
The above experiments show that, 
to have a small QMC-error by low WAFOM point sets, 
the integrand should have high order partial
derivatives with small norms
(see a preceding research \cite{HARASE-OHORI}, too).
However, WAFOM seems to work with some non-differentiable
functions (such as Continuous in the previous section).

\subsection{$t$-value again}
Niederreiter-Pirsic \cite{NIEDERREITER-PIRSIC}
showed that for a digital net $P$,
the strict $t$-value of $P$ as a $(t,m,s)$-net is expressed 
as  
\begin{equation}\label{eq:t-value}
m-t+1 = \min_{A \in P^\perp-\{0\}}\mu_1(A).
\end{equation}
Here $\mu_1$ is 
Dick's $\alpha$-weight
for $\alpha=1$, which is known as
the Niederreiter-Rosenbloom-Tsfasman weight.

There is a strong resemblance between (\ref{eq:t-value})
and Definition~\ref{def:WAFOM}.
Again in (\ref{eq:t-value}), high complexity of all elements in $P^\perp-\{0\}$
gives strong uniformity (i.e., small $t$-value).
The right hand side of (\ref{eq:t-value}) 
is efficiently computable by a MacWilliams-type
identity in $O(sN\log N)$ steps of integer operation
\cite{DICK-MATSUMOTO}.
\begin{question}
The formula (\ref{eq:t-value}) for $t$-value uses
the minimum over $P$, while Definition~\ref{def:WAFOM}
of WAFOM and (\ref{eq:W-alpha}) use the summation over $P$.
Can we connect 
$t$-value in (\ref{eq:t-value})
with WAFOM in Definition~\ref{def:WAFOM}?
It may perhaps relate with ultra-discretization \cite{ULTRA-DISCRETIZATION}.
\end{question}

\section{Randomization by digital shift}
Let $P \subset M_{s,n}(\Ftwo)$ be a linear subspace.
Choose $\sigma \in M_{s,n}(\Ftwo)$.
The point set 
$P+\sigma :=\{B+\sigma | B \in P \}$
is called the digital shift of $P$ by $\sigma$.
Since $P+\sigma$ is not an $\Ftwo$-linear subspace,
one can not define $\wafom(P+\sigma)$. 
Nevertheless, the same error bound holds as $P$.
Under a uniform random choice of $\sigma$, $P+\sigma$
becomes unbiased. 
Moreover, the mean square error is bounded as follows: 
\begin{theorem} (Goda-Ohori-Suzuki-Yoshiki \cite{GOSY})
$$
\Error(f_n;P+\sigma) \leq C(s,n)||f||_n \wafom(P), \mbox{ and}
$$
$$
\sqrt{{\mathbb E}(\Error(f_n;P+\sigma)^2)} \leq 
C(s,n)||f||_n \wafom^{\mbox{r.m.s.}}(P),
$$
$$ \mbox{where }
\wafom^{\mbox{r.m.s.}}(P):=
\sqrt{\sum_{A \in P^\perp-\{0\}}2^{-{2}\mu(A)}}.
$$
\end{theorem}

\section{Variants of WAFOM}
As mentioned in the previous section, \cite{GOSY} defined
$\wafom^{\mbox{r.m.s.}}(P)$. 
As another direction, the following generalization 
of WAFOM is proposed by Yoshiki \cite{YOSHIKIBOUND}
and Ohori \cite{OHORI-MASTER}: in Definition~\ref{def:Dick-weight}, 
the function $\mu(A)$ might be generalized by:
$$
\mu_\delta(A):=\sum_{1\leq i \leq s, 1\leq j \leq n} (j+\delta)a_{ij}
$$
for any (even negative) real number $\delta$ (note that
this definition is different from that of $\mu_\alpha$, 
but we could not find a better notation).
Then Definition~\ref{def:WAFOM} gives $\wafom_\delta(P)$.
The case where $\delta=1$ is dealt in \cite{YOSHIKIBOUND}.
A weak point of the original WAFOM 
is that WAFOM value does not vary enough 
and consequently it is not useful in 
grading point sets for a large $s$, see Figure~2, the $s=8$ case.
By choosing a suitable $\delta$, we obtain $\wafom_\delta(P)$ that 
varies for large $s$ (even for $s=16$) and useful in choosing a good
point set
\cite{OHORI-MASTER}.
A table of bases of such point sets is 
available from Ohori's GitHub Pages: 
\url{http://majiang.github.io/qmc/index.html}.
These point sets are obtained by Ohori, using 
Harase's method based on linear scrambling, from NX sequences.
Thus, they have small $t$-values and small WAFOM values.
Experiments show their good performance \cite{MORI-MULTI}. 

\section{Conclusion}
Walsh figure of merit (WAFOM) \cite{WAFOM} for $\Ftwo$-linear point sets
as a quality measure for a QMC rule is discussed.
Since WAFOM satisfies a Koksma-Hlawka type inequality
(\ref{eq:boundbyWAFOM}), its effectiveness for very smooth functions
is assured. Through the experiments on QMC integration, we observed that 
the low WAFOM point sets show higher order convergence
such as $O(N^{-1.2})$ for several test functions (including non-smooth one)
in dimension four, and $O(N^{-1.05})$ for dimension eight.

\begin{acknowledgement}
The authors are deeply indebted to Josef Dick, who
patiently and generously informed us of beautiful researches
in this area, and to Harald Niederreiter for leading us
to this research. They thank for the indispensable helps
by the members
of Komaba-Applied-Algebra Seminar (KAPALS): Takashi Goda, Shin Harase, 
Shinsuke Mori, Syoiti Ninomiya, Mutsuo Saito, Kosuke Suzuki, and Takehito Yoshiki.
We are thankful to the referees, who informed of numerous improvements
on the manuscript.
The first author is partially supported by 
JSPS/MEXT Grant-in-Aid for Scientific Research 
No.21654017, No.23244002, No.24654019, and No.15K13460.
The second author is partially supported by the Program for Leading
Graduate Schools, MEXT, Japan.
\end{acknowledgement}

%
\bibliographystyle{spmpsci}
\bibliography{sfmt-kanren}
%
%
\end{document}